\documentclass[reqno]{amsart}

\usepackage{mathbbol,mathrsfs}
\usepackage{amsfonts,amssymb,amsmath,amsthm,amscd}
\usepackage{latexsym}
\usepackage[mathscr]{eucal}
\usepackage{wrapfig,graphics,graphicx,epic,epsfig}
\usepackage{ifthen}
\usepackage{oldgerm,units}
\usepackage{pstricks,pst-node}

 \input xy
\xyoption{all}

 \makeatletter
\renewcommand*\env@matrix[1][*\c@MaxMatrixCols c]{%
  \hskip -\arraycolsep
  \let\@ifnextchar\new@ifnextchar
  \array{#1}}
\makeatother

\newcommand{\joverline}[2]{%
  \mathord{% make sure we're in math mode
    \vbox{\offinterlineskip
      \halign{##\cr
        $\scriptscriptstyle#1$\hrulefill\cr
        \noalign{\kern.4ex}
        $\; #2 \,$\cr
      }%
    }%
  }%
}
\newtheorem{theorem}{Theorem}[section]

\newtheorem{example}[theorem]{Example}

\newcommand\GLnA[1]{{\operatorname{GL}}{[1]}}
\newcommand\SLA[1]{{\operatorname{SL}}{[1]}}

\def\semirings0{semirings$^\dagger$}
\def\Semirings0{Semirings$^\dagger$}
\def\semiring0{semiring$^\dagger$}
\def\hsemiring0{$\frac 12$-semiring$^\dagger$ with negation}
\def\nsemiring0{semiring$^\dagger$ with negation}
\def\nSemirings0{Semirings$^\dagger$ with negation}
\def\hsemirings0{$\frac 12$-semirings$^\dagger$ with negation}
\def\hSemirings0{$\frac 12$-Semirings$^\dagger$ with negation}
\def\domain0{domain$^\dagger$}

\def\domains0{domains$^\dagger$}
\def\field0{semifield$^\dagger$}
\def\fields0{semifields$^\dagger$}

\newcommand{\etype}[1]{\renewcommand{\labelenumi}{(#1{enumi})}}
\def\eroman{\etype{\roman}}

\def\({\left(}
\def\){\right)}

\def\inv{{^{-1}}}

\newcommand{\trop}[1]{\mathcal{#1}}

\newcommand{\tT}{\trop{T}}

\def\one{\mathbb{1}}
\def\zero{\mathbb {0}}

\def\rone{\one_R}
\def\rzero{\zero_R}

\def\R{\mathbb R}

\newtheorem{thm}[theorem]{Theorem}
\newtheorem*{thm*}{Theorem}
\newtheorem*{dig*}{Digression}
\newtheorem{lem}[theorem]{Lemma}
\newtheorem{rem}[theorem]{Remark}
\newtheorem{prop*}{Proposition}

\newtheorem{prop}[theorem]{Proposition}
\newtheorem{defn}[theorem]{Definition}
\newtheorem*{examp*}{Example}
\newtheorem*{examples*}{Examples}
\newtheorem*{remark*}{Remark}

\newtheorem*{defn*}{Definition}
\newtheorem*{note*}{Note}
\DeclareMathAlphabet{\mathbbold}{U}{bbold}{m}{n}
\newcommand{\mzero}{\mathbbold{0_M}}

\begin{document}
\title[Hypergroups and hyperfields]{Hypergroups and hyperfields in universal algebra}

%******************************* AMS classification ***********************
%\subjclass[2000]{Primary 15A09, 15A03, 15A15, 65F15; Secondary
%16Y60 }

%******************************* date *************************************
\author[L.~Rowen]{Louis Rowen}
\address{Department of Mathematics, Bar-Ilan University, Ramat-Gan 52900,
Israel} \email{rowen@math.biu.ac.il}

%******************************* AMS classification ***********************
\subjclass[2010]{Primary    16Y60,     12K10, 06F05, 14T05
Secondary 12K10, .}

%******************************* date *************************************
\date{\today}

%******************************* keywords *********************************

\keywords{ hyperfield, negation map,   tropical algebra, tropical
geometry, power set, semiring, supertropical
algebra.} %% section %%%%%%%%%%%%%%%%%%%%%%%%%%%%%

%%%%%%%%%%%%%%%%%%%%%%%%%%%%%%%%% section %%%%%%%%%%%%%%%%%%%%%%%%%%%%%

%{\small \tableofcontents}
\numberwithin{equation}{section}
%%%%%%%%%%%%%%%%%%%%%%%%%%%%%%%% section %%%%%%%%%%%%%%%%%%%%%%%%%%%%%
\begin{abstract} Hypergroups are lifted to power semigroups with negation, yielding
a method of transferring results from semigroup theory. This applies
to analogous structures such as hypergroups,  hyperfields, and
hypermodules, and permits us to transfer the general theory espoused
in \cite{Ro2} to the hypertheory.
\end{abstract}
\maketitle

\section{Introduction}

This note, a companion to \cite{Ro2}, grew out of a conversation
with Matt Baker, in which we realized that the ``tropical
hyperfield'' of \cite{Bak} and \cite[\S5.2]{Vi} is isomorphic to the
``extended'' tropical arithmetic in Izhakian's Ph.D. dissertation
(Tel-Aviv University) in 2005, also
cf.~\cite{zur05TropicalAlgebra,IzhakianRowen2007SuperTropical}. On
the other hand, there are many parallels between the two theories.
This motivated us to see whether hyperfields in general also can be
studied by semiring theory, which fits in well with the theory of
universal algebra, and which might be more amenable for further
study. Viewing a hyperfield as a group with additive structure on
part of its power set, we want to extend this structure to all of
the power set, thereby making the definitions tighter and
``improving'' the additive structure to make standard tools more
available for developing an algebraic theory.

Thus, the theme is to embed the category of hyperfields (and their
modules) into the category of semirings with negation (and their
modules), as studied in \cite{Ro2}, defined on power sets. The
tricky part in passing to the power set is distributivity, which
must be weakened at times to a notion that we call ``weak
distributivity,'' and we thereby weaken ``semiring'' to
``$\tT$-semiring,'' where distributivity holds only with respect to
a special subset $\tT$. Then we can
treat all hyperrings (not just hyperfields) in this context. %In
%particular, given a hyperfield $R$ we pass to the sub-semiring
%$\tilde R$ of the power set ``$\tT$-semiring,'' which is called a
%``partial semiring'' in \cite{Ro2}.

It turns out that the hyperrings of \cite{Bak,Vi} can be injected
naturally into $\tT$-semirings, which are power sets with a negation
map, in the context of \cite{Ro2}, whereby the hyperring is
identified with the subset of singletons. Then one can develop
linear algebra over hyperfields, and also go through \cite{Bak},
making the appropriate adjustments to view matroids over these
semifields with negation.

In the other direction, Henry~\cite{He} has defined a hypergroup
structure on symmetrized monoids.

\subsection{Pre-semirings}$ $

Semigroups need not have   an  identity element, but monoids do, and
are usually written multiplicatively, i.e.~$(A, \cdot , \one)$. We
start with semirings, for which the standard reference is
\cite{golan92}. Since distributivity will be weakened, we remove it
(as well as the element $\zero$) from the definition. We enrich a
monoid
 with a ``hyper'' operation (resembling addition).

\begin{defn} A monoid $(A,\cdot,\one)$ \textbf{acts} on a set $\mathcal S$
if there is a multiplication $A \times \mathcal S \to \mathcal S$
satisfying $\one s = s$ and $(a_1a_2)s = a_1 (a_2 s)$ for all $a_i
\in A$ and $s \in S.$

 A \textbf{pre-semiring} $(A,\cdot,+,\one)$ is a multiplicative
monoid $(A, \cdot, \rone )$ also possessing the structure of an
additive Abelian semigroup $(A, +)$, on which $(A,\cdot)$ acts.

A \textbf{pre-semifield}   is a pre-semiring $(A,\cdot,+,\one)$ for
which $(A,\cdot, \one)$ is a group.

 A \textbf{premodule} $\mathcal S$ over a monoid $(A,\cdot,\one)$ is an
Abelian group $(\mathcal S, +, \zero)$ on which $A$ acts, also
satisfying the condition:

 $$r\mzero = \mzero, \ \forall r \in R.$$

 A \textbf{premodule} $\mathcal S$ over a pre-semiring $(A,\cdot,+,\one)$ is a
premodule over the monoid $(A,\cdot,\one)$, also satisfying the
condition:

\begin{itemize}
  \item ``Left distributivity'':  $
   (a_1+a_2)s = a_1s + a_2s,$ $\forall a_i \in A,$ $\forall s \in S.$
\end{itemize}
\end{defn}

Then a \semiring0 $(R, +, \cdot, \rone )$ (without $\zero) $ would
be a pre-semiring satisfying  the usual distributive laws.  There
are various versions of distributivity that will be relevant to us
later.

\begin{defn} Define the following notions, where $\mathcal (S,+)$ is a
premodule over a pre-semiring $(A,\cdot,+,\one)$, for $a_i \in
A$ and $s_j\in S$:
\begin{enumerate}\eroman
  \item ``Right distributivity'':
$   a(s_1+s_2) = as_1 + as_2.$
   \item  ``Two-sided distributivity'': Left and right distributivity.
   \item ``Double distributivity'':  $
   (a_1+a_2)(s_1+s_2) = a_1s_1 + a_2s_1 + a_1s_2 + a_2 s_2.$
  \item ``Generalized distributivity'':  $$\left(\sum _i a_i\right)\left(\sum _j s_j\right) = \sum_{i,j}
 (a_is_j).$$
\end{enumerate}
\end{defn}

\begin{rem} The   properties are in  order of increasing formal strength although, by
induction, double distributivity implies generalized distributivity.
%
%Distributivity in a
% semiring implies generalized distributivity. But we will encounter
% structures for which there is a difference.
\end{rem}

\subsection{Motivation: Power sets of semigroups}\label{powermon}$ $

We let $\mathcal P (A)$ denote the  power set of a set $A$, i.e.,
the set of subsets of $A$. The sets $\{ a \}$ for $a\in A$ are
called \textbf{singletons}.

\begin{thm}\label{mon1} $ $ \begin{enumerate}\eroman
   \item
     Given a
 monoid $(A,\cdot,\one)$, we can extend its operation
  to $\mathcal P(A)$ elementwise, by putting
$$S_1  S_2 =\{{  s_1  s_2 :\ s_j \in S_j} \}.$$
Then $(\mathcal P(A),\cdot)$ also is a monoid with the identity
element $\{ \one\},$ on which $A$ acts.

   Given a
 semigroup $(A,+)$, we can define addition
 elementwise on $\mathcal P(A)$ by defining
$$S_1 + S_2 =\{{  s_1 + s_2 :\ s_j \in S_j}  \}.$$
Then $(\mathcal P(A),+)$ also is a semigroup.

Thus, when $A$ is a pre-semiring (resp.~semiring), so is $\mathcal
P(A)$, and $\mathcal P(A)$ is an $A$-premodule via the action $$aS =
\{ as : s \in S\}.$$

 \item More generally, the relevant concepts in  universal algebra were
 outlined in \cite[\S 2.3]{Ro2}, including signatures defined via operators and identical relations. We can lift operators from an $(\Omega; \operatorname{Id})$-algebra
 $\mathcal A$
 to $\mathcal P(\mathcal A),$ as follows:
Given an operator $\omega_{m } = \omega_{m }(x_{1 }, \dots, x_{m })$
on $\mathcal A$, we define $\omega_{m }$ on $\mathcal P(\mathcal A)$
via $$ \omega_{m }(S_{1 }, \dots, S_{m }) = \{ \omega_{m }(s_{1 },
\dots, s_{m }) :  s_{k }\in  S_{k }, \, 1 \le k \le m\}.$$

Let us call an identical  relation
 $\phi(x_1,\dots, x_\ell)  =  \phi(x_1,\dots, x_\ell)$ \textbf{multilinear}
 if each $x_i$ appears exactly once in the definition (via the operators).
  Then any multilinear identical relation holding in
$\mathcal A$ also holds in $\mathcal P(\mathcal
A)$.\end{enumerate}\end{thm}
\begin{proof}
(i) First we verify associativity:
\begin{equation}\begin{aligned} (S_1+S_2) + S_3  & =  \{a_1 + a_2: a_j \in S_j\}
+ S_3 \\& = \{(a_1 + a_2) + a_3 : a_j \in S_j\} \\& = \{a_1 +(a_2 +
a_3) : a_j \in S_j\} \\& =  S_1 + \{a_2 + a_3: a_j \in S_j\}  \\& =
 S_1+(S_2 + S_3).
\end{aligned}\end{equation}

For generalized distributivity we have:

\begin{equation} \begin{aligned}  \sum S_i \sum T_j & =  \left\{\sum _{i_k} \sum _{k=1}^{m_i} a_{i_k}: a_{i_k} \in S_i\right\}
 \left\{\sum _{j_\ell} \sum _{\ell=1}^{n_j} b_{j_\ell}: b_{j_\ell} \in T_j\right\} \\ &  =   \left\{ \sum _{i,j,k,\ell} a_{i_k}   b_{j_\ell}: a_{i_k}
 \in S_i,\
b_{j_\ell} \in T_j\right\}  =
 \sum S_i T_j    .
 \end{aligned}\end{equation}

(ii) We generalize the proof in (i). By an easy induction applied to
\cite[Definition~2.12]{Ro2}, any formula $\phi(x_1,\dots, x_\ell)$
satisfies
$$\phi(S_1,\dots, S_\ell)  = \{ \phi(s_1,\dots, s_\ell) : s_j \in
S_j\},$$ and thus any multilinear identical relation $\phi = \psi$
holding elementwise in $\mathcal A$ also holds set-wise in $\mathcal
P(\mathcal A)$.
 \end{proof}

In particular, distributivity and generalized distributivity lift
from a semiring $A$ to $\mathcal P(A)$. When the relation is not
multilinear we encounter difficulties due to repetition. For
example, being a group does not lift/ Indeed, the defining identical
relation $x x^{-1} = 1$ is quadratic in $x$; here we are defining
the inverse as a unary operation $\omega_1: x \mapsto x^{-1},$ so
the left side is $x \omega _1(x).$ In fact the only invertible
elements in $\mathcal P(A)$  are the singletons.

\subsection{Hypermonoids  and semirings}$ $

 The next step is to formulate all of our extra structure in terms
of addition (and possibly other operations) on $\mathcal P(A)$ as a
premodule over    $A$. But this is not as easy
when $A$ itself is not a semiring, so let us pause to review
hypermonoids, to see just how much of the semiring structure we
would need.

The ``intuitive'' definition: A  hypermonoid should be a triple
$(A,\boxplus,\zero)$ where $\boxplus : A \times A \to \mathcal
P(A)$, and the analog of associativity holds:
$$ (a_1 \boxplus a_2) \boxplus a_3 = a_1 \boxplus (a_2 \boxplus
a_3), \quad \forall a \in A. $$

There is a fundamental difficulty in this definition ---  $\quad a_1
\boxplus a_2$ is a set, not an element of $A$, so technically $ (a_1
\boxplus a_2) \boxplus a_3 $ is not defined. This difficulty is
exacerbated when considering generalized associativity; for example,
what does $ (a_1 \boxplus a_2) \boxplus (a_3\boxplus a_4) $ mean? We
rectify this by passing to $\mathcal P(A).$

\begin{defn} A \textbf{hypermonoid} is a triple
$(A,\boxplus,\zero)$ where
\begin{enumerate}\eroman
   \item
$\boxplus$ is a commutative binary operation $A \times A \to
\mathcal P(A),$ which also is associative in the sense that if we
define $$a \boxplus S = \cup _{s \in S} \ a \boxplus s,$$ then $(a_1
\boxplus a_2) \boxplus a_3 = a_1 \boxplus (a_2\boxplus a_3)$ for all
$a_i$ in $A.$
  \item
$\zero$ is the neutral element.
\end{enumerate}

We always think of $\boxplus$ as a sort of addition.\end{defn}

 We write
$\tilde A$ for $\{ a_1 \boxplus a_2 : a_i \in A\}.$ Note that $\{ a
\} = a \boxplus \zero \in \tilde A.$ Thus there is a natural
embedding $A \hookrightarrow \tilde A$ given by $a \mapsto \{ a \},$
and we can transfer the addition to $\mathcal P(A)$ by defining
$$ \{ a_1 \} \boxplus \{ a_2 \} =  a_1 \boxplus  a_2.$$

By definition, the hypermonoid is not closed under repeated
addition, which makes it difficult to check basic identical
relations such as associativity.

 Many hypermonoids satisfy the extra   property:

\bigskip

\textbf{Property P}. $a,b \in a\boxplus b$ whenever $a\boxplus b$ is
not a singleton.

\bigskip

A \textbf{hyperinverse} of an element $a$ in a hypermonoid
$(A,\boxplus,\zero)$ is an element denoted as $-a$, for which $\zero
\in a \boxplus (-a).$

 A \textbf{hyperzero} of a hypermonoid $(A,\boxplus,\zero)$ is an element of the form $a \boxplus (-a) \subseteq \mathcal P(A).$

A \textbf{hypergroup} is a hypermonoid $(A,\boxplus,\zero)$,
satisfying the extra property:

 \bigskip

  \noindent Every element $ a \in A$ has a unique hyperinverse.

 \bigskip

 By \cite[\S 2]{He}, any hypergroup   satisfies the  condition:

 \bigskip

 \noindent (\textbf{Reversibility})\label{rever} If $ a \in b \boxplus
 c,
\quad \text{then} \quad c \in a \boxplus (-b).$

 \bigskip In \cite[Definition 3.1]{Vi}
 Viro~
 calls this a \textbf{multigroup}.

\begin{defn} A \textbf{hypermodule} over a monoid $(A,\cdot,\one)$ is a
hypermonoid $(M,\boxplus,\zero)$ together with an action of $A$ on
$M$ such that
   distributivity holds for $M$ over $A$.

 A   \textbf{hyperring} $(
 A,\cdot, \boxplus,\zero)$ is a hypermonoid  $(
 A,\boxplus,\zero)$ which also is a hypermodule over  $(A,\cdot,\one)$.
    \end{defn}

In other words,
 a
 hyperring  is  an additive hypermonoid  which also is a monoid
 with respect to an associative multiplication that distributes over
addition; we have the two operations $\cdot$ on $A$ and $\boxplus: A
\to \tilde A,$ with distributivity holding on the elements of $A$.
%Note that we cannot yet discuss double distributivity or generalized
%distributivity at the level of the hypermodule or hyperring since we
%do not have the sum of two sets.

 \bigskip

A \textbf{hyperfield} is a hyperring $(
 A,\cdot, \boxplus,\zero)$, with $(A, \cdot)$
a group.

 \bigskip

 \cite[Definition 2.3]{Ja} defines a \textbf{hypermonoid morphism}
 to be a map $f: A_1 \to A_2$ of hypergroups, satisfying
 $f( a
\boxplus  b) \subseteq f( a) \boxplus  f(b).$ This yields the
category  of hypergroups and their morphisms, which matches the
definition of morphism in \cite{Ro2}.
%
%\begin{lem}{mon11}
%   Given a
% hypergroup $(A,+)$, we can define addition
% elementwise on $\mathcal P(A)$ by defining
%$$S_1 \boxplus S_2 =\bigcup {{  s_1 \boxplus s_2 :\ s_j \in S_j} }.$$
%Then $(\mathcal P(A),+)$ also is a semigroup. If $A$ has an identity
%element $0,$ then $(\mathcal P(A),+)$ has the identity element $\{
%0\}.$
%\end{lem}
%\begin{proof} As in Theorem~\ref{mon1}, since we are given that
%associativity holds at the level of elements.
%\end{proof}

 \bigskip

Here are some easy instances in which associativity fails in
$\mathcal P(A)$.

\begin{example}
 Consider the natural max-plus algebra.

\begin{enumerate}\eroman
   \item Define $a\boxplus b = \sup \{a,b\}$,
$a \ne b,$ and $$a\boxplus a = \{ 0 , 9\}.$$ Then each element has a
unique hyperinverse, itself, and this is associative on distinct
single elements (taking their max) but $(2\boxplus 2) \boxplus  5 =
\{ 0, 9\} \boxplus  5 = \{ 5,9\}  $ whereas $2\boxplus (2 \boxplus
5) = 5.$
 \item Define $a\boxplus b = \sup \{a,b\}$,
$a \ne b,$ and $$a\boxplus a = \{-a, 0 ,a\}.$$   Again each element
has a unique hyperinverse, itself, and now
$$(2\boxplus 3) \boxplus 3 = 3 \boxplus 3 = \{ -3, 0 , 3 \};$$
$$2\boxplus (3 \boxplus 3 ) = 2 \boxplus \{ -3, 0 , 3 \} =  \{ 2, 2 , 3
\}.$$
\end{enumerate}
\end{example}

 \begin{lem} If $(A,\cdot,\one)$ is a monoid and $(A,\boxplus,\zero)$
 is a hypermonoid, then $A$ acts on $\tilde A$  via the action
\begin{equation}\label{prem} aS = \{ as : s \in S\}.\end{equation}.
\end{lem}

\begin{proof} $(a_1a_2)S = \{ (a_1a_2) s: s\in S\} =   \{ a_1(a_2s): s\in
S\} = a_1(a_2S).$

 \end{proof}

\subsection{The power set of a hyperfield}$ $

To proceed further, we need associativity at the level of sets, and
we need the following definition to make this precise (and hopefully
more manageable, since then we can do all the calculations in the
power set).

\begin{defn} Every operator $\omega$ on $A$ is extended to an \textbf{element-compatible} operator on $\mathcal P(A)$, in the sense that
 \begin{equation}\label{elcom}\omega_{m}(S_{1}, \dots,
S_{m}) = \bigcup \{ \omega_{m}(s_{1}, \dots, s_{m}) : s_{k}\in
S_{k}, \, 1 \le k \le m\}.\end{equation}
\end{defn}

If $(A,\boxplus,0)$ is a hypermonoid, then $\mathcal P(A)$ is  a
monoid with respect to a commutative associative binary operation
$\mathcal P(A) \times \mathcal P(A) \to \mathcal P(A),$  compatible
with the  operation on singletons, in the sense that
\begin{equation}\label{elcom1} S_1 \boxplus S_2 = \bigcup \left\{\{
s_1\} \boxplus \{ s_2\} : s_i \in S_i\right\}\end{equation}  $\tilde
A$ is the subset of $\mathcal P(A)$ containing all singletons and
their sums. ($\{\zero\}$ is the neutral element.)

\begin{lem}\label{mon51} For any multiplicative monoid $A$, any  invertible element of $\mathcal P(A)$
must be a singleton.
\end{lem}
\begin{proof} If $S$ has two elements $s_1,s_2$ and is invertible,
then $S^{-1}$ contains $s_1\inv$ and $s_2\inv$, implying $1 =
s_1\inv s_2,$ i.e., $s_1=s_2$.
 \end{proof}

Theorem~\ref{mon1} generalizes easily to:

\begin{thm}\label{mon2}    Given a
 hyperfield $(A,\boxplus,\zero )$, we can define addition
 elementwise on $\mathcal P(A)$ by means of \eqref{elcom1}.
Then $(\mathcal P(A),\boxplus)$ is a monoid, whose identity element
is $\{ \zero  \}.$  In this case $A \setminus \{ \zero  \} $ (viewed
as the set of singletons) is the set of invertible elements of
$\mathcal P(A)$.
\end{thm}
\begin{proof} We need to verify associativity, repeating the proof
of Theorem~\ref{mon1}, replacing $+$ by $\boxplus$.
\begin{equation}\begin{aligned} (S_1\boxplus S_2) \boxplus  S_3  & =  \bigcup \{s_1 \boxplus  s_2: s_j \in S_j \}
\boxplus  S_3 \\& = \bigcup ((s_1 \boxplus  s_2) \boxplus  s_3 ):
s_j \in S_j\} \\& = \bigcup (s_1 \boxplus (s_2 \boxplus  s_3)) : s_j
\in S_j\}
\\& = S_1 \boxplus  \bigcup \{s_2 \boxplus  s_3: s_j \in S_j\}
\\& =
 S_1\boxplus (S_2 \boxplus  S_3).
\end{aligned}\end{equation}

 The set of invertible elements of $\mathcal P(A)$ must be contained in the set of singletons of $\mathcal
 P(A)$, which is  $A$.

 \end{proof}

 For any finite set $S = \{ s_1, \dots, s_m\}\subset \mathcal P(A),$ we write
 $\boxplus S$ for $s _1\boxplus \dots \boxplus s_m,$ which makes
 sense since we already have   associativity of $\boxplus$.
%
%  Given a
% multiplication $(A,\cdot)$, we  define multiplication
% elementwise on $\mathcal P(A)$ by defining
%$$S_1  S_2 =\bigcup {\{  s_1  s_2 :\ s_j \in S_j\} };$$

\begin{rem}\label{tri}
 Viro showed that the general transition of universal
relations to $\mathcal P(A)$ is not as straightforward as it may
seem. The analogous argument to Theorem~\ref{mon2} unravels for
hyperrings, since distributivity does not pass from elements to
sets:

\begin{enumerate}\eroman
   \item\cite[Theorem~4.B]{Vi} $(a \boxplus b)(c \boxplus d) \subseteq (ac)
\boxplus (ad) \boxplus (bc) \boxplus (bd)$ in any hyperring; the
same argument shows that $(\boxplus S)(\boxplus T) \subseteq
\boxplus (ST)$ for any finite sets $S,T$;

 \item \cite[Theorem~5.B]{Vi} Recall Viro's  ``triangle'' hyperfield,
  defined over $\R ^+$ by the
 formula $$a \boxplus b = \{c \in \R^+ :
|a - b| \le c \le a + b\}.$$ In other words, $c \in  a \boxplus b$
iff there exists an Euclidean triangle with sides of lengths $a, b,$
and $ c.$

The ``triangle'' hyperfield $R$ does not
 satisfy ``double distributivity,'' so $\mathcal P(R)$ is not
 distributive.
 \end{enumerate}
\end{rem}

 In other words, the analog of Theorem~\ref{mon1} fails for distributivity. To overcome
 this setback, we need to modify our underlying algebraic structure
 both
 at the hyper level and the power set level,
 the crux of the matter being to weaken generalized distributivity on
 sets.

Actually, many of the important examples of hyperfields are doubly
distributive, so we could pass to (distributive) power semirings
without further ado. But even in the absence of doubly
distributivity, we can formulate the weaker theory in terms of
universal algebra, in order to have those techniques at our
disposal. On the face of it, this is problematic since the hypersum
set could be arbitrarily large. However, we can get around this by
focusing on the   monoid of singletons, and using operators instead
of elements.

\subsection{Pre-semirings}

Motivated by the fact that a hyperring is a multiplicative monoid,
we bring in the following definition, keeping the power set in mind:

\begin{defn}\label{modu2} A \textbf{pre-\semiring0} is a set $(A, +, \cdot, \rone
)$ for which $(A, +)$ is  an additive Abelian semigroup and $(A,
\cdot, \rone )$ is a multiplicative monoid but not necessarily
satisfying the usual distributive laws.

Given a  pre-\semiring0 $A$ and a (distinguished) multiplicative
submonoid $\tT$, an $(A,\tT)$-\textbf{module}  is a premodule over
$A$ that also satisfies the distributivity conditions for $r_i \in
\tT$ and $a_i \in M$:

 \begin{enumerate}\eroman
 \item $(r_1 + r_2) a =  r_1a + r_2 a, $
  \item $r ( a_1 + a_2) =  r a_1 + r a_2, $
\end{enumerate}
\end{defn}
 In line
with \cite{IzhakianRowen2007SuperTropical}, we call $\tT$
\textbf{tangible}. When these conditions hold for $A$, we then call
$(A,\tT)$ a \textbf{$\tT$-\semiring0}.
%
%Define $\widehat{\tT} = \{ \boxplus S: S \subseteq T \}.$
%$\widehat{\tT}$ is a $ \tT $-premodule under the given
%multiplication $\tT \times \widehat{\tT} \to \widehat{\tT}$, but not
%necessarily a \semiring0 because of lack of generalized
%distributivity.

This can be described in the framework of universal algebra, with
elaboration and details given in \cite[\S~5]{Ro2}. We define the
left multiplication maps $\ell_r : M \to M$ by $\ell _r(a) = ra,$ a
unary operator for each $r \in \tT,$ and rewrite these rules as the
identical relations

 \begin{enumerate}\eroman
 \item $\ell_{r_1 + r_2}(x) = \ell_{r_1}(x) + \ell_{r_2}(x), $
  \item $\ell_r ( x_1 + x_2) =  \ell_r (x_1) + \ell_r (x_2) $
\end{enumerate}

\begin{example}
For any hyperring $\tT$, taking  $P(\tT), \tT)$ is a
$\tT$-\semiring0.
  \end{example}

\subsubsection{The categorical approach}$ $

We can improve these results slightly by working in a category with
a weaker definition of morphism.

 \begin{defn}  Multiplication \textbf{weakly distributes} over
addition in a subset $\mathcal S \subseteq \mathcal P(\tT)$ if for
all finite $S,T \subseteq \mathcal S$   we have
\begin{equation}\label{weakdist0} (\boxplus S)(\boxplus T)   \subseteq \boxplus (S T).\end{equation}
In this case we call $\mathcal P(\tT)$ a \textbf{weak power
semiring}. (It is not a semiring.)
\end{defn}

%
%%
% \begin{defn} Suppose $(A, \cdot)$ is a multiplicative monoid   and
% $(\mathcal P(A),\boxplus,
% \{ 0 \}$ is a power monoid. Multiplication \textbf{weakly distributes} over
%addition in a subset $\mathcal S \subseteq \mathcal P(A)$ if for all
%finite $S,T \subseteq \mathcal S$   we have
%\begin{equation}\label{weakdist0} (\boxplus S)(\boxplus T)   \subseteq \boxplus (S T).\end{equation}
%\end{defn}

The restriction of this definition to hyperrings  is:

 \begin{defn} Suppose $(\tT, \cdot)$ is a monoid and
  $(\tT, \boxplus)$ a hypermonoid.
  Multiplication \textbf{weakly distributes} over
$\boxplus$ in $\tT$ if for all  $a_i,b \in  \tT$   we have
\begin{equation}\label{weakdist} (\boxplus _i a_i)b   \subseteq \boxplus _i (a_ib).\end{equation}
In this case we call $\tT$ a \textbf{weak hyperring}.

%In this case, when $(\mathcal S, \boxplus, \{ 0 \})$ is closed under
%$\boxplus$ and multiplication, we call $\mathcal S$ a \textbf{weak
%sub-semiring} of $\tT$.
\end{defn}

\begin{prop}\label{mon3} Suppose $\tT$ is a hyperring. Then multiplication  weakly distributes  over
$\boxplus$ in $                \mathcal P(\tT)$.
\end{prop}
\begin{proof} We need to verify \eqref{weakdist0}.
But writing $S = \{ s_1, \dots, s_m\}$  we have $$(\boxplus
S)(\boxplus T) = \bigcup _i (\boxplus s_i) T =  \bigcup _i \boxplus
(s_i T) \subseteq \boxplus (S T)
$$ since each $s_i T \subseteq  \boxplus ST.$
 \end{proof}

The reverse inclusion fails since we simultaneously encounter
varying $s_i T$ when $i$ varies. Thus, we are interested in weak
power semirings which, strictly speaking,   are not quite semirings.
So far, the power set of a hyperring is a weak power semiring. Now
we repeat the proof of \cite[Theorem~4.B]{Vi}, to show that at the
bottom level we have not lost anything.

\begin{thm}\label{invertit} Suppose  $\tT$ is a  weak hyperring. If $(\tT,\cdot)$ also is a
 group, then
   $\tT$ is a hyperring, which can be identified with the set of singletons of $\mathcal P(\tT)$.
\end{thm}
\begin{proof} As in \cite[Theorem~4.A]{Vi}, to obtain distributivity, we need to reverse the inequality
\eqref{weakdist0} when $S = \{a \}$ is a singleton $\{ a\}$, given
multiplicative inverses in~$\tT$.  But  $$ \boxplus (aT) =  a
a^{-1}(\boxplus (aT) )  \subseteq a (\boxplus a^{-1} aT) = a
\boxplus  T  .$$ %The last assertion now is immediate from
%\cite[Theorem~4.B]{Vi}, since the product of singletons is a
%singleton by hypothesis.
\end{proof}

%
%Technically, the restriction of a weak power semiring $(\mathcal
%P(A),\boxplus,
% \{ 0 \}$
% to $\{\a : a \in A\}$
% is a weak hyperring, but not necessarily a hyperring. But
%when $A$ is a group we get a hyperfield by \cite[Theorem~4.A]{Vi}.
Thus, the theory of hyperfields  embeds into the theory of weak
power semirings.

\begin{defn} A  \textbf{weak morphism} of weak power semirings
 is a multiplicative homomorphism $f: \mathcal P(\tT_1) \to \mathcal P(\tT_2)$,
 satisfying
  $$f( S_1 \boxplus S_2
  ) \subseteq f(S_1) \boxplus  f(S_2).$$

%This yields the category  $\wsr$  of weak power semirings and their
%morphisms.
\end{defn}

By induction, we have $f(\boxplus  S)\subseteq  \boxplus f(S)$ for
any finite set $S$, and thus for arbitrary $S$. This is described in
universal algebra in \cite[Definition~5.4]{Ro2}. Note that the
multiplicative version yields equality, since if $f(\{a\}) f(\{ b\})
 \subseteq
f(\{ab\}),$ then $f(\{a\} )f(\{b\})
 =f(\{ab\})$ since they are
both singletons, so $f( a  )f( b )
 =f( ab )$.

\subsection{Weak power modules}$ $

As often is the case, one gets a deeper understanding by turning to
modules. This was done in \cite[Definition~2.18]{Bak}, but we take a
slightly weaker definition in line with our earlier categorical
considerations.

\begin{defn}\label{modu1}
Suppose that $\mathcal P(\tT)$ is a weak power semiring. A
 \textbf{weak module}  over $\mathcal P(\tT)$ is a set $M$ with
 an element $\mzero$ such that  $(\mathcal P(M),\boxplus_M,\mzero)$ is
 a   monoid with scalar multiplication
$ \tT  \times \mathcal P(M)\to \mathcal P(M)$ satisfying the
following for all $S _i \subseteq \tT$ and $T \subseteq M$:
\begin{enumerate}\eroman
   \item $(S_1S_2)T = S_1(S_2T)$;
  \item $(\boxplus _\tT S)( \boxplus _M T )\subseteq \boxplus (ST)$;
 \item  $r\mzero = \mzero$
for all $r$ in $\tT$;
   \item $\rzero a = \mzero$, $\forall
a \in M.$\end{enumerate}
\end{defn}

This leads to a slight modification for hypermodules:

\begin{defn}\label{modu2}
 A \textbf{weak hypermodule} over a weak hyperring $(R,\boxplus_R,\rzero)$ is
a hypergroup $(M,\boxplus_M,\zero)$ together with a binary operation
$R \times M \to \mathcal P(M)$ satisfying the following properties
for all $r, r _i \in R$ and $a, a_j \in M$:
\begin{enumerate}\eroman
   \item $(r_1r_2)a = r_1(r_2a)$;
  \item $(\boxplus_i r_i) (\boxplus_j a_j) \subseteq \boxplus_{i,j}
  (r_ia_j);$
 \item  $r\mzero = \mzero$;
   \item $\rzero a = \mzero$.\end{enumerate}
\end{defn}

 As before,  $\{ \boxplus S: \text{ finite }S
\subseteq M\}$ is a submodule of $\mathcal P(M)$.

\begin{defn} A  \textbf{weak module morphism}
 is a map $f: \mathcal P(M) \to \mathcal P(N)$ of weak modules,
 satisfying $$f( r S_1)
\subseteq r f( S_1);$$
 $$f( S_1
\boxplus  S_2) \subseteq f( S_1) \boxplus  f(S_2), \quad \forall S_i
\subseteq M.$$
\end{defn}

The weak module morphism is an example of a morphism as given in
\cite{Ro2}.

\begin{thm}\label{invertit2} Suppose $\mathcal P(R)$ is a weak power semiring
and  $(R,\cdot)$ also is a
 group. If $M$  is a weak module  over $\mathcal P(R)$, then
  $M$ is an
$R$-module.
\end{thm}
\begin{proof} We repeat the proof of Theorem~\ref{invertit}. First we show that $f( r S) = r f( S),$
i.e., $r f(s) = f(rs) \in f(rS).$ for all $r\in R$ and $s\in S$.
Indeed,
$$ r f(s) = r f(r\inv rs) \subseteq r r\inv f(rs) = f(rs).$$

 To obtain distributivity, we need to reverse the inclusion~\eqref{weakdist0},
given multiplicative inverses in $R$. Taking $S = \{a \}$ to be a
singleton, we are given $$ \boxplus (aT) =  a a^{-1}(\boxplus (aT))
\subseteq a (\boxplus (a^{-1} aT)) = a \boxplus  T  .$$
\end{proof}

Taking $f$ to be left multiplication by an element
 $a\in R$, we have  $a (\boxplus  S) \subseteq  \boxplus
(aS)$ which is precisely weak distributivity. In other words, $M$ is
weakly distributive over $R$ iff left multiplication is a  weak
module morphism, for every element $r\in R$. When $(R,\cdot)$ is a
group, $M$ is  distributive over $R$.

\subsubsection{Negation maps}$ $

We also want to treat negation maps from \cite[\S 4]{Ro2} in this
perspective. We review the definition.

\begin{defn}\label{negmap}
 A \textbf{negation map} on an additive semigroup $(A,+)$ is a
 semigroup homomorphism
$(-) :A \to A$ of order $\le 2,$ again  written $a\mapsto
(-)a$.\end{defn}

(Thus $(-)(a+b) = (-)a + (-)b.)$ For all other operators, including
multiplication, we have a different perspective:

\begin{defn}\label{neg0}  A \textbf{negation map} $(-) :
\mathcal A \to \mathcal A$ on an operator $\omega_{m,j}$ (other than
addition)
  satisfies
\begin{equation}\label{neg1} (-)\omega_{m,j}( a_1, \dots,  a_{u-1}, a_u, a_{u+1},\dots,  a_m)=
 \omega_{m,j}(a_1, \dots, a_{u-1}, (-)a_u, a_{u+1},
\dots, a_m)\end{equation} for each  $a_u \in \mathcal A_u,$ $1\le
u\le m$.

% In particular, a \textbf{negation map} on a {\it multiplicative} monoid
%$(A,\cdot)$ is a map $(-):A \to A$ of order~$\le 2,$ written
%$a\mapsto (-)a,$ satisfying $(-)((-)a)=a$ for all $a \in A,$ as well
%as
%\begin{equation}\label{neg} (-)(a_1a_2)  = ((-)a_1) a_2 = a_1((-)a_2).\end{equation}
%
%A \textbf{negation map} on a semiring $(R,\cdot,+,\rone)$ is a map
%$(-):R\to R$  which is simultaneously a negation map on the additive
%semigroup $(R,+)$ and on the  multiplicative monoid $(R,\cdot)$
\end{defn}

\begin{lem}\label{mon5}  Any negation map on $\mathcal A $
induces a negation map on $\mathcal P(\mathcal A)$, via $(-)S = \{
(-)s: s \in S \}.$  \end{lem}
\begin{proof} A special case of Theorem~\ref{mon1}(ii), viewing the
properties of the negation map as identical relations.
 To see that $(-)(a_1+a_2) =  (-)a_1 (-) a_2,$
note that $ 0 \in a_i (-)a_i$ for $i=1,2,$ so $$0 \in   a_1 (-)a_1
+a_2 (-) a_2 = (a_1+a_2) (-)a_1 (-) a_2.$$

Likewise, $ 0 \in a_1 (-)a_1$ implies $ 0 \in a(a_1 (-)a_1)  \in
aa_1 (-)aa_1.$
%
%In case  $A$ is a hyperring we note from the previous paragraph that
%$(-)(a_1a_2)= ((-)a_1)a_2 = a_1((-)a_2),$ and thus $((-)a_1)((-)a_2)
%= (-) ((-)a_1) a_2 = a_1 a_2.$
 \end{proof}

%
%\begin{example}
% Consider the natural max-plus algebra. Define $a+b$ as usual unless
%$a=b,$ and $$a+a = \{ 0 , 9\}.$$ Then each element has a unique
%inverse, itself, so this seems to be a hypergroup which is
%associative on distinct single elements (taking their max), but
%$(2+2) + 5 = \{ 0, 9\} + 5 = \{ 5,9\}$ whereas $2+(2 + 5) = 5$

%\end{example}

When $\tT$ is a hyperfield,  $(A,\tT)$ is a
  $\tT$-\semiring0  with the  negation
  map   $a \mapsto (-)a.$
In this way  we have embedded the theory of hyperfields into the
theory of $\tT$-\semirings0 with a negation map. This is pushed even
further in
  \cite[\S7.9]{Ro2}.

%
%
%\subsubsection{
%Summary of the main notions}$ $
%
%Given a hypermodule (resp.~hyperring) $M$, we want to work in the
%weak power module $\mathcal P(M)$ with its corresopnding negation
%map (resp.~ weak power semiring $\mathcal P(R)$ with negation map,
%which is an example of a $\tT$-\semiring0). $M$ is identified with
%the singletons of $\mathcal P(R)$, which we called the
%\textbf{tangible} elements; these apply for example to the tropical
%theory. In order to handle associativity of addition we need to pass
%to the distinguished additive monoid $\tilde M = \{\boxplus S: |S| <
%\infty \}.$
%
%We are given a negation map $(-)$ on $\mathcal P(M)$, and define
%$$\mathcal P(M) ^{\circ}= \{ a \boxplus ((-)a): a \in M\},$$
%the $M$-submodule of quasizeros of $\mathcal P(M).$ (When $M$ is a
%hyperring then $\mathcal P(M) ^{\circ}$ is an ideal.) This will take
%the place of zero when extending results from classical algebra.
%
%We also have the \textbf{surpassing} relation $s_1 \preceq S_2$ if
%$s_1 \boxplus s^\circ \subseteq S_2$ for some $s \in \mathcal P(M),$
%and $S_1 \preceq S_2$ if $s_1 \in S_2,$ $\forall s_1 \in S_1,$
%treated abstractly in \cite{Ro2}.

\subsection{Major examples}$ $

Let us see how all of this applies to the major examples of
\cite{Bak}. Since these examples are so important, we will pay
special attention to the set   $\tilde \tT$ corresponding to a
hyperring $\tT$. Although the theory presented above formally passes
to the weak power semiring $\mathcal P (\tT)$, one actually gets
distributivity in $\mathcal P (\tT)$ when the underlying hyperring
satisfies generalized distributivity, which happens in many of the
examples. Even better, we can often identify $\mathcal P(\tT)$ with
a semiring which we already recognize.

Many of the ``good'' examples can be put in the framework of
\cite[Remark 2.7]{Bak}.

\begin{example}\label{Basicexample0} Let $R $
be a commutative semiring.  Any multiplicative monoid $\tT$,
together with a surjection of multiplicative monoids $\varphi : R
\to \tT $, has an induced hyperring structure given by the
hyperaddition law
$$a_1 \boxplus  a_2 := \varphi (\varphi ^{-1} (a_1) + \varphi ^{-1}
(a_2)).$$ This extends naturally to $\mathcal P(\tT)$, via
$$S_1
\boxplus  S_2 := \varphi (\varphi ^{-1} (S_1) + \varphi ^{-1}
(S_2)).$$ Generalized distributivity on $\mathcal P(\tT)$ and thus
on $\tilde \tT$, is inherited from generalized distributivity on
$\mathcal P(R).$ Explicitly, for $  a_i \in  S$ and $  b_j \in T,$
we have
\begin{equation}\label{dist22} (\boxplus _i a_i) (\boxplus _j b_j) =\sum _i \varphi (\varphi ^{-1} (a_i))\sum _j \varphi (\varphi ^{-1}
(b_j)) = \sum_{i,j} \varphi (\varphi ^{-1} (a_i)\varphi ^{-1} (b_j))
= \sum_{i,j} \varphi (\varphi ^{-1} (a_i b_j))\subseteq \boxplus (S
T),\end{equation} yielding  $\boxplus S \boxplus T \subseteq
\boxplus (S T)$. For the opposite direction, given $\boxplus _{i,j}
a_i b_j \in \boxplus (ST),$ we reverse \eqref{dist22} to get
$$\sum_{i,j} \varphi (\varphi ^{-1} (a_i b_j))=\sum _i \varphi (\varphi ^{-1} (a_i))\sum _j \varphi (\varphi ^{-1}
(b_j))\in (\boxplus S)(\boxplus T).$$
 Thus $\mathcal P(\tT)$ is a semiring,
and its theory can be embedded into semiring theory.
\end{example}

The complications arise when Example~\ref{Basicexample0} is not
applicable, cf.~(vii), (viii) of the next example.

\begin{example}\label{Basicexamples}$ $

 The tropical hyperfield.
Define  $\R _\infty = \R \cup \{ - \infty\}$ and define the  product
$a \bigodot b := a + b$ and
$$a \boxplus b=
\begin{cases} max(a, b)\text{ if } a \ne b, \\   \{ c : c \le a\}
\text{ if } a = b.
\end{cases}$$
 Thus 0 is the multiplicative identity element, $-
\infty$ is the additive identity, and we have a hyperfield
(satisfying Property P), called the \textbf{tropical hyperfield}.

\begin{prop} This is easily seen to be isomorphic (as hyperfields) to
Izhakian's \textbf{extended tropical arithmetic}
\cite{zur05TropicalAlgebra},
  further expounded as \textbf{supertropical algebra} in
\cite{IzhakianRowen2007SuperTropical}, where we identify
$(-\infty,a] : =\{ c : c \le a\}$ with $a^\nu$, so   we have a
natural hyperfield isomorphism of this tropical hyperfield with the
sub-semiring $\widehat{\R _\infty}   $ of $\mathcal P(\R _\infty)$,
because

$$(-\infty,a] + b = \begin{cases} b : b >a;\\ (-\infty,a] : b =a \\  (-\infty,b] \cup (b,a] = (-\infty,a] : b
<a.
\end{cases} $$

\end{prop}

This isomorphism is   as semirings. Thus Example~\ref{Basicexample0}
is applicable.
\end{example}

\begin{example}
  The Krasner hyperfield. Let $K = \{ 0; 1 \}$  with the usual
operations of Boolean algebra, except that now $ 1 \boxplus  1 = \{
0; 1 \} .$ The Krasner hyperfield satisfies Property P. Again, this
generates a sub-semiring of $\mathcal P (K)$ having three elements,
and is just the supertropical algebra of the monoid~$K$, where we
identify $\{ 0; 1 \} $ with $ 1^\nu.$  Example~\ref{Basicexample0}
is
 applicable.
\end{example}

\begin{example} Valuative hyperfields (\cite[Example 2.12]{Bak}) also are
isomorphic to the extended semirings in the sense of
\cite{IzhakianRowen2007SuperTropical}, in the same way.

\end{example}

\begin{example} (Hyperfield of signs) Let $ S := \{ 0, 1 , -1\}$  with the usual
multiplication law and hyperaddition defined by $1 \boxplus  1 = \{
1\} ,$ $-1 \boxplus  -1 = \{ -1\} ,$ $ x \boxplus  0 = 0 \boxplus  x
= \{ x\} ,$ and $1 \boxplus  -1 = -1 \boxplus  1 = \{ 0, 1,-1\} .$
Then $S$ is a hyperfield (satisfying Property P), called the
hyperfield of signs.

In this case, though, we have a natural interpretation for $\mathcal
P(S):$

\begin{enumerate}\eroman
  \item $\{+1\}$ means ``positive,'' denoted as $>_0$.

 \item $\{ -1\}$ means ``negative,'' denoted as $<_0$.

 \item $\{ 0\}$ means ``neutral.''

\item $\{0,+1\}$ means ``non-negative,'' denoted as $\ge _0$.

\item $\{0,-1\}$ means ``non-positive,'' denoted as $\le _0$.

\item $S = \{-1,0,1\}$ means ``could be anything''

(Note that we have not denoted $\{-1, +1\}$.)
\end{enumerate}

Then we have the familiar identifications:

\begin{enumerate}\eroman
  \item  $>_0 + 0 \ = \ >_0+  \ge _0 \ = \ >_0+>_0 \ = \ >_0;$

  \item  $\ge _0 + 0 \ = \ \ge _0+  \ge _0 \ = \  \ge _0;$

 \item  $0 + 0 \ = \ 0;$

 \item  $\le _0 + 0 \ = \ \le _0+  \le _0 \ = \  \le _0;$

 \item  $<_0 + 0 \ = \ <_0+  \le _0 \ = \ <_0+<_0 \ = \ <_0;$

  \item $>_0+<_0 \ =\ \ge_0+<_0 \ =\ >_0+\le_0 \ =\
  \ge_0+\le_0 \ =\ S + <_0 \ = \ S + 0 \ = \ S + >_0 \ =  \ S + \le_0 \  \ =  \ S + \ge_0 \ = \ S.$

 These six elements constitute the sub-semiring $\tilde S$ of $\mathcal
P(S)$.

\end{enumerate}
\end{example}

\begin{example}The ``triangle'' hyperfield  $A$ of Remark~\ref{tri} is not doubly distributive but does satisfy
Property P since $|a-b| \le a \le a+b.$ Here $\tilde A  = \{
[a_1,a_2]: a_1 \le a_2 \}$ since $[a_1,a_2] = \frac {a_1+a_2}2 +
\frac {a_2-a_1}2 \in \tilde A $.  Any interval $[0,b]$ is in
$\mathcal P(A)^\circ,$ since $[0,b] = [0,\frac b 2] + [0,\frac b
2].$

\end{example}

\begin{example} The phase hyperfield. Let $S^1$ denote the complex
unit circle, and $ P := S^1 \cup \{ 0 \} .$ We say that points $a$
and $b$ are antipodes if $a = -b.$ Multiplication is defined as
usual (so corresponds on~$S^1$ to addition of angles). We call an
arc of less than 180 degrees \textbf{short}. The hypersum is given
by
$$a \boxplus b=
\begin{cases} \text{ all points in the short arc from }a  \text{ to } b  \text{ if } a \ne b;\\   \{ -a,0,a \} \text{
if } a , b \text{ are antipodes;     } \\   \{  a \} \text{ if }   b
= 0 .
\end{cases}$$  Then $P$ is a hyperfield (satisfying Property P), called the \textbf{phase
hyperfield}. At the power set level, given $T_1,T_2 \subseteq S^1$,
one of which having at least two points, we define $T_1 \boxplus
T_2$ to be the union of all (short) arcs from a point of $T_1$ to a
non-antipodal point in $T_2$ (which together makes a connected arc),
together with  $ \{ 0 \}$ if $T_2$ contains an antipode of $T_1$.
Note that any arc of $S^1$ can be obtained by taking $T_1$ to be a
single point in the middle and $T_2$ to be the two endpoints.

In other words, $\widetilde{S^1}$ of $\mathcal P(S^1)$  is the set
of short arcs, possibly with $\{ 0 \}$ adjoined. This is not a
sub-semiringl for this we need the set of all arcs, where $\boxplus$
is concatenation (and filling in the rest of $S^1$ if the arcs go
more than half way around), and adjoining $\{ 0 \}$ if the arcs
contain an antipode.

Double distributivity fails, when we take $a_1$ and $a_2$ almost to
be antipodes, $b_1 = a_2,$ and the arc connecting $b_1$ and $b_2$
just passes the antipode of  $a_1$; then  $(a_1 \boxplus a_2)(b_1
\boxplus b_2)$ is the arc from $a_1$ to $b_2,$ a little more than a
semicircle, whereas $a_1 b_1  \boxplus a_1 b_2  \boxplus a_2 b_2$ is
already all of $S^1.$

 Viro~\cite{Vi} also has a somewhat different version,
 which is not distributive.
\end{example}

\begin{example} Here is another example, suggested by Lopez, also cf.~\cite{Lop}. Consider $\R$,
with addition given by $a \boxplus b $ and $b \boxplus a $ (for $a
\le b$) to be the interval $[ a,b]$. This extends to  addition on
intervals, given by $[ a_1,b_2] + [ a_2,b_2] = \{ \min (a_1,a_2),
\max (b_1,b_2)\},$ which clearly is  associative. But the inverse is
not unique, since $a + (-a) = [-a,a]$ contains 0, but so does $\frac
a 2 + a.$ On the other hand, this does satisfy the restriction that
every set of the form $a + (-a)$ cannot be of the form $a + (-b)$
for $b \ne a,$ so if we modify the condition of quasi-inverse to
stipulate that $a + (-a)$ must be of the form $c+ (-c)$ for some
$c$, then it is unique. This is essentially the general condition
set forth in \cite{Ro2}.

\end{example}


\begin{thebibliography}{10} %{IMS}




\bibitem{AGG}
M.~Akian,   S.~Gaubert, and A.~Guterman.
\newblock Linear independence over tropical semirings and beyond.
\newblock In {\em Tropical and  Idempotent Mathematics}, G.L. Litvinov and S.N.
Sergeev, (eds.),
\newblock {\em Contemp. Math.},  495:1--38, 2009.
%\bibitem{BR} \newblock  Computational aspects of polynomial identities  \newblock
%A.K.~Pet

  \bibitem{AGG1}
M.~Akian,   S.~Gaubert, and A.~Guterman. Tropical Cramer
determinants
    revisited. In Tropical and idempotent mathematics and
    applications,
   Contemp. Math. 616, 1--45. Amer. Math. Soc., Providence, RI,
   2014.

\bibitem{Bak} M.~Baker. Matroids over hyperfields,  arXiv:1601.01204v2 [math.CO], 2016.

\bibitem{BlR} G.~Blachar and L.~Rowen {Symmetrized rank of matrices}

\bibitem{CC} A.~Connes and C.~Consani, From monoids to
hyperstructures: in search of an absolute arithmetic. In {\em
Casimir Force, Casimir operators, and the Riemann hypothesis},,
Walter de Gruyter, Berlin, 147--198, 2010.



%\bibitem{Ei} R.~Eilhauer, \newblock {\em Zur Theorie der Halbk\"{o}rper I, Acta Math. Sci. Acad. Hung 19 (1968), 23--45.

\bibitem{G.Ths}
  S.~Gaubert, {\em Th\'{e}orie des syst\`{e}mes lin\'{e}aires
dans les dio\"{i}des}. PhD dissertation, School of Mines. Paris,
July  1992.

\bibitem{golan92}
J.~Golan, \newblock {\em Semirings and their Applications},
Springer-Science + Business, Dordrecht, 1999.
\newblock (Previously published by Kluwer Acad. Publ.,  1999.)

  \bibitem{GoM}     M.~Gondran and M.~Minoux. Graphs, dioids and semirings, volume 41 of Operations
    Research/Computer Science Interfaces Series. Springer, New York, 2008.

\bibitem{He} Henry

\bibitem{zur05TropicalAlgebra}
Z.~Izhakian.
\newblock Tropical arithmetic and matrix algebra.
\newblock {\em Comm. in Algebra}  {37}(4):1445--1468, {2009}.

\bibitem{IzhakianRowen2007SuperTropical}
Z.~Izhakian and L.~Rowen.
\newblock {Supertropical algebra}.
\newblock   {\em Adv. in Math.}, 225(4):2222--2286, 2010,
 Preprint at arXiv:0806.1175, 2007.
%
%%6
%\bibitem{IKR-LinAlg2}  Z.~Izhakian, M.~Knebusch, and   L.~Rowen.
%\newblock Dual spaces and bilinear forms in supertropical linear
%algebra, {\em Linear and Multilinear Algebra}, 60(7):865--883,
%2012.
%%. (Preprint  at arXiv:1201.6481, 2011.)
%
%%7
%\bibitem{IzhakianRowen2009TropicalRank}
%Z.~Izhakian and L.~Rowen.
%\newblock The tropical rank of a tropical matrix.
%\newblock {\em Comm. in Algebra},  {37}(11):{3912--3927}, 2009.
%%
%%%9
%\bibitem{IzhakianRowen2008Matrices}
%Z.~Izhakian and L.~Rowen.   \newblock  {Supertropical matrix
%algebra.} \newblock  {\em Israel J. of Math.},
%  182(1):383--424, 2011.

%\bibitem{IzhakianRowen2009Equations}
%Z.~Izhakian and L.~Rowen.
% Supertropical matrix algebra II: Solving tropical
%equations, \emph{Israel J. of Math.}, 186(1): 69--97, 2011.
% % (Preprint at arXiv:0902.10.302159, 16pp.)

%
%\bibitem{IzhakianRowen2010MatrixIII}
%Z.~Izhakian and L.~Rowen. \newblock
% Supertropical matrix algebra III: Powers
%of matrices and their supertropical eigenvalues, \emph{J.~ of
%Algebra}, 341(1):125--149, 2011.
% %  (Preprint at arXiv:1008.0023.)

\bibitem{Ja} J.~Jaiung, \newblock {Algebraic geometry over hyperrings}.
Preprint. Available at arxiv:math.AG/1512.04837, 37 pages, 2015.

\bibitem{Lop} M.~G\'omez,
S.~L\'opez-Permouth , F.~Mazariegos, A.J.~ Vargas De Le\'on, and
R.Z.~Cifuentes, Group Structures on Families of Subsets of a Group,
    arXiv:1604.01119 [math.GR], 2015.



\bibitem{Niv} A.~Niv. \newblock {On pseudo-inverses of matrices and their
characteristic polynomials in supertropical algebra}, \newblock
\emph{Linear Algebra Appl.}  471: 264--290, 2015.



\bibitem{Niv2}
A. Niv. \newblock {Factorization of tropical matrices}, \newblock
\emph{J.~Algebra Appl.}, 13(1):1350066:1--26, 2014.

\bibitem{Pl} M. Plus. Linear systems in (max; +)-algebra. In Proceedings of the 29th Conference on
    Decision and Control, Honolulu, Dec. 1990.


\bibitem{RS} C.~Reutenauer and H.~Straubing, \textit{Inversion of matrices over a
commutative semiring}, J.~Algebra, \textbf{88}, 350–-360,  1984.

\bibitem{Ro}
L.H.~Rowen. {\em Graduate algebra: A noncommutative view}
\newblock {\em Semigroups and Combinatorial Applications}.
\newblock American Mathematical Society, 2008.

\bibitem{Ro2}
L.H.~Rowen. {\em Algebras with a negation map},  originally entitled
{\em Symmetries in tropical algebra}, arXiv:1602.00353 [math.RA],
2016.

\bibitem{Rut} D.E.~Rutherford. Inverses of Boolean matrices,
    Proceedings of the Glasgow Mathematical Association 6 (1963),  49--53.

    \bibitem{St} H.~Straubing. A combinatorial proof of the Cayley-Hamilton
theorem. Discrete Math., 43(2-3):273–-279, 1983.

\bibitem{Vi} O.Y.~Viro. Hyperfields for Tropcial Geometry I.
Hyperfields and dequantization. ArXiv AG/1006.3034, 2010.

 \end{thebibliography}
\end{document}